\documentclass{amsart}

\usepackage{graphicx}

\newtheorem{theorem}{Theorem}
\newtheorem{corollary}[theorem]{Corollary}
\newtheorem{proposition}[theorem]{Proposition}
\newtheorem{claim}{Claim}

\theoremstyle{remark}
\newtheorem{remark}{Remark}

\newcommand{\mindeg}{\mathop{\mathrm{mindeg}}\nolimits}

\begin{document}

\title{Cyclic and finite surgeries on Montesinos knots}

\author{Kazuhiro Ichihara}
\address{School of Mathematics Education,
Nara University of Education,
Takabatake-cho, Nara, 630-8528, Japan}
\email{ichihara@nara-edu.ac.jp}
\urladdr{http://mailsrv.nara-edu.ac.jp/~ichihara/index.html}
\thanks{The first author is partially supported by
Grant-in-Aid for Young Scientists (B), No. 20740039,
Ministry of Education, Culture, Sports, Science and Technology, Japan.}

\author{In Dae JONG}
\address{Graduate School of Science, Osaka City University, Osaka 558-8585, Japan}
\email{jong@sci.osaka-cu.ac.jp}
\urladdr{http://www.ex.media.osaka-cu.ac.jp/~d07sa009/index.html}

\date{\today}


\begin{abstract}
We give a complete classification of 
the Dehn surgeries on Montesinos knots 
which yield manifolds with cyclic or finite fundamental groups. 
\end{abstract}
\maketitle

\section{Introduction}

A \textit{Dehn surgery} on a knot $K$ in a $3$-manifold $M$ 
is an operation to create a new $3$-manifold 
from $M$ and $K$ as follows: 
Remove an open tubular neighborhood of $K$, 
and glue a solid torus back. 
By gluing a solid torus back as it was, 
the surgery gives the original manifold again. 
So such a surgery is called \textit{trivial}, and 
we will ignore it in general.

On knots in the $3$-sphere $S^3$, 
it is an interesting problem to determine and classify 
all non-trivial Dehn surgeries which produce 
$3$-manifolds with cyclic or finite fundamental groups, 
which we call \textit{cyclic surgeries} / \textit{finite surgeries} respectively.

As part of the Hyperbolic Dehn Surgery Theorem, 
Thurston \cite{Thurston} established that 
there are finitely many cyclic and finite surgeries. 
In fact, Culler, Gordon, Luecke and Shalen \cite{CGLS} 
(respectively, Boyer and Zhang \cite{BoyerZhang}) proved 
there are at most three cyclic (resp., five finite) surgeries. 
Furthermore, it is conjectured that knots admitting 
cyclic (resp., finite) surgeries are doubly primitive 
(resp., primitive/Seifert fibered) knots as introduced by 
Berge \cite{Berge} (resp., Dean \cite{Dean}). 
See \cite[Problem 1.77]{Kirby} for more information.

Cyclic and finite surgeries have been studied extensively 
for some classes of knots. 
For example, 
it was shown by Delman and Roberts in \cite{DelmanRoberts} that  
no hyperbolic alternating knot admits a cyclic or finite surgery. 

\begin{figure}[htb]
\begin{center}
 \begin{picture}(100,80)
  \put(0,0){\scalebox{0.2}{\includegraphics{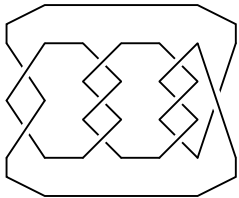}}}
 \end{picture}
 \caption{A diagram of a Montesinos knot}
 \label{fig:MontesinosKnot}
\end{center}
\end{figure}

One of the other well-known classes of knots, 
containing non-alternating ones, is the Montesinos knots. 
A \textit{Montesinos knot} is defined as a knot 
admitting a diagram obtained
by putting rational tangles together in a circle.
See Figure \ref{fig:MontesinosKnot} for instance.
In particular, a Montesinos knot $K$ is called 
a \textit{$(a_1,a_2,\cdots,a_n)$-pretzel knot} if 
the rational tangles in $K$ are of the form $1/{a_1},1/{a_2},\cdots,1/{a_n}$. 

In this paper, 
based on studies by Delman \cite{Delman2} and Mattman \cite{Mattman}, 
we give a complete classification of 
cyclic / finite surgeries on Montesinos knots as follows. 

\begin{theorem}\label{MainThm}
Let $K$ be a hyperbolic Montesinos knot. 
If $K$ admits a non-trivial cyclic surgery, 
then $K$ must be equivalent to 
the $(-2,3,7)$-pretzel knot and the surgery slope is $18$ or $19$. 
If $K$ admits a non-trivial acyclic finite surgery, 
then $K$ must be equivalent to either 
the $(-2,3,7)$-pretzel knot and the surgery slope is $17$, or
the $(-2,3,9)$-pretzel knot and the surgery slope is $22$ or $23$. 
\end{theorem}

As a direct corollary, together with the result by Wu \cite{Wu1}, 
we have the following. 

\begin{corollary}\label{cor}
Let $K$ be a hyperbolic arborescent knot. 
If $K$ admits a non-trivial cyclic surgery, 
then $K$ must be equivalent to 
the $(-2,3,7)$-pretzel knot and the surgery slope is $18$ or $19$. 
If $K$ admits a non-trivial acyclic finite surgery, 
then $K$ must be equivalent to either 
the $(-2,3,7)$-pretzel knot and the surgery slope is $17$, or
the $(-2,3,9)$-pretzel knot and the surgery slope is $22$ or $23$. \qed
\end{corollary}

Recently, using Khovanov homology, 
it was shown in \cite[Theorem 7.5]{Watson} that 
the $(-2, p, p)$-pretzel knot does not admit finite surgeries 
for $p \in \{ 5, 7, \cdots , 25 \}$. 

Very recently, Futer, Ishikawa, Kabaya, Mattman, and Shimokawa 
\cite{fiKMS} obtained, independently, 
a complete classification of finite surgeries 
on $(-2, p, q)$-pretzel knots with odd positive integers $p$ and $q$.

\begin{remark}
It is already known 
which Montesinos knots are non-hyperbolic. 
If a Montesinos knot is equivalent to 
one consisting of at most two rational tangles, 
then it actually is a two-bridge knot. 
Menasco \cite{Menasco} showed that 
the non-hyperbolic two-bridge knots are the $(2, p)$-torus knots. 
The only other non-trivial non-hyperbolic Montesinos knots are 
the $(-2, 3, 3)$- and $(-2, 3, 5)$-pretzel knots, 
which are actually the $(3, 4)$- and $(3, 5)$-torus knots, respectively. 
This was originally shown by Oertel \cite[Corollary 5]{Oertel} 
as well as in an unpublished monograph \cite{BonahonSiebenmann} 
by Bonahon and Siebenmann. 
The cyclic and finite surgeries of 
torus knots have been completely classified by Moser \cite{Moser}.  
\end{remark}

To prove Theorem~\ref{MainThm}, 
we will prepare two propositions, 
Propositions \ref{prop1} and \ref{prop2}, 
which will be shown in Sections \ref{sec2} and \ref{sec3} respectively. 
Then, in the last section, 
Theorem \ref{MainThm} will be proved from these propositions 
together with a result of Mattman \cite{Mattman}.

\subsection*{Acknowledgments}
The authors would like to thank Mikami Hirasawa 
for useful conversations about the fiberedness of pretzel knots 
as in the proof of Claim~\ref{clm:fibered}. 
They also wish to thank Hiroshi Matsuda 
for his useful comments 
about the rank and the Euler characteristic 
of Heegaard Floer homology 
as in the proof of Claim~\ref{clm:Lsurgery}, 
and also thank Thomas Mattman 
for his useful comments about 
cyclic and finite surgeries on arborescent knots 
as in Corollary~\ref{cor}. 
They also thank the referee for careful reading of our draft 
and useful suggestions to improve the proof of Claim 2.

\section{Cyclic/finite surgeries and the Alexander polynomials}\label{sec2}

In this section, we prove the following proposition. 

\begin{proposition}\label{prop1}
Let $K$ be a hyperbolic Montesinos knot 
admitting a non-trivial cyclic or finite surgery. 
Then $K$ is equivalent to a $(-1, 2n, p, q)$-pretzel knot, 
where $n$ is a non-zero integer and 
$p,q$ are odd positive integers with $3 \le p \le q$. 
Furthermore 
all non-zero coefficients of the Alexander polynomial 
for $K$ are $\pm 1$. 
\end{proposition}

\begin{proof}
Suppose that a hyperbolic Montesinos knot $K$ 
admits a non-trivial cyclic or finite surgery. 
Then Delman showed in \cite{Delman, Delman2} that 
$K$ must be equivalent to either 
a $(-2l, p, q)$-pretzel knot, 
a $(-1, 2n, p, q)$-pretzel knot or a $(-1, -1, 2m, p, q)$-pretzel knot 
with an integer $n$, integers $l, m>1$ 
and odd positive integers $p,q$ ($3 \le p \le q$). 
Also see \cite[Section 2, Section 3]{Wu2}. 
Actually Delman showed that 
any Montesinos knot except for those listed above
admits an essential lamination in its exterior 
which survives all non-trivial Dehn surgeries. 
\textit{Essential laminations} were introduced by Gabai and Oertel 
in \cite{GabaiOertel} and, actually, they showed that if a 3-manifold admits 
an essential lamination, then its universal cover must be the 
3-space $\mathbb{R}^3$. 
In particular its fundamental group is not cyclic or finite. 
See \cite{GabaiOertel} for the precise definition. 

By virtue of Delman's result, 
in order to prove Proposition \ref{prop1}, 
it suffices to show that 
the first and the third types of pretzel knots described above cannot have cyclic or finite surgeries. 
Note here that 
a $(-2, p, q)$-pretzel knot (the case $l=1$ in the first) 
is equivalent to 
a $(-1, 2, p, q)$-pretzel knot (the case $n=1$ in the second). 
Also a $(-1, -1, 2, p, q)$-pretzel knot 
(the case $m=1$ in the third) 
is equivalent to 
a $(-1, -2, p, q)$-pretzel knot (the case $n=-1$ in the second). 
Thus, excluding over laps, we are assuming $l, m \ne 1$. 

Among the classes of knots described above, 
the first one was already studied by Mattman in \cite{Mattman}. 
He actually showed in \cite[Theorem 1.1 and 1.2]{Mattman} that 
any $(-2l, p, q)$-pretzel knot with $l >1$ 
and odd positive integers $p,q$ ($3 \le p \le q$) 
has neither cyclic surgeries nor finite surgeries. 

Thus, in the following, we focus on the third class above. 
We here use the following strong result by Ni, \cite[Corollary 1.3]{YiNi},  
established by using the Heegaard Floer homology theory: 
If a knot in $S^3$ admits a cyclic or finite surgery, 
then it must be a fibered knot.
Actually he showed that a knot $K$ in $S^3$ must be fibered 
if $K$ admits a surgery yielding an L-space. 
Here a rational homology sphere $Y$ is called an \textit{L-space} 
if the rank of $\widehat{HF}(Y )$ is equal to $| H_1 (Y ; \mathbb{Z}) |$. 
In fact, any $3$-manifold with a cyclic or finite fundamental group 
is an L-space, as is shown in \cite[Proposition 2.3]{OzsvathSzabo2}. 

Now, the next claim, together with the result by Ni, imply 
the first conclusion of Proposition \ref{prop1}. 

\begin{claim}\label{clm:fibered}
Let $m > 1$ and $p, q$ be odd positive integers $(p \le q)$. 
The $(-1, -1, 2m, p, q)$-pretzel knot is not fibered.
\end{claim}

\begin{proof}
We just apply the algorithm given in \cite[Theorem 6.7]{Gabai}.
Here we include only an outline, 
assuming that the reader is rather familiar with \cite[Theorem 6.7]{Gabai}. 
Please see \cite{Gabai} for details. 

Let $K$ be a $(-1,-1, 2m, p, q)$-pretzel knot 
with an integer $m>1$ and odd positive integers $p,q$ ($p \le q$). 
We start to apply the algorithm in \cite[Theorem 6.7]{Gabai} 
with $n_1 = -1$, $n_2 = -1$, $n_3 = 2m$, $n_4 = p$, $n_5=q$. 
After a cyclic permutation, 
the surface $R$ obtained by applying 
Seifert's algorithm is of type II in \cite[TYPE II.6.5]{Gabai} 
with $m_1 = -1$, $m_{11}=2m$, $m_2 = p$, $m_3 = q$, $m_4 = -1$. 
(See \cite[Figure 6.3]{Gabai}.) 
We now see CASE 2 in \cite[Theorem 6.7]{Gabai}. 
Here we note that 
the associated oriented pretzel link $L'$ 
(defined in \cite[TYPE II.6.5]{Gabai}) 
is of type $( 2m , -2 ,-2) $. 
Since $\displaystyle \sum_{j=1}^4 \frac{m_j}{|m_j|} = -1+1+1-1 = 0$ and 
$L'$ is of type $( 2m , -2 ,-2) \ne \pm (2,-2,2)$ if $m>1$, 
we check CASE 2B in \cite[Theorem 6.7]{Gabai}. 
Then we see that $K$ is fibered if and only if $L'$ is fibered. 
For $L'$, we check CASE 1 in \cite[Theorem 6.7]{Gabai}, 
and verify that $L'$ is not fibered since 
no $n_j$ is $\pm1$ and $L'$ is not equivalent to 
a pretzel link of type $\pm(2,-2,\cdots,2,-2,n)$ with an integer $n$. 
Therefore we conclude that $K$ is not fibered. 
\end{proof}

The second conclusion of Proposition \ref{prop1} 
follows from results of Ozsv\'{a}th and Szab\'{o}, 
also achieved by using the Heegaard Floer homology theory. 
We first prepare the following claim, 
which is implicitly used in \cite[Proof of Corollary 1.3]{YiNi}. 

\begin{claim}\label{clm:Lsurgery}
If $\alpha / \beta$-Dehn surgery on a non-trivial knot $K$ in $S^3$ yields an L-space 
for some coprime integers $\alpha , \beta$ with $\beta \ge 2$, 
then $\alpha$-Dehn surgery on $K$ also yields an L-space. 
\end{claim}

\begin{proof}
Given coprime integers $\alpha , \beta$ and a knot $K$ in $S^3$, 
let $S^3_{\alpha / \beta} (K)$ denote the 3-manifold 
obtained from $S^3$ by $\alpha / \beta$-surgery on $K$. 
We recall the following general formula 
given in \cite[Proposition 9.5]{OzsvathSzabo3}:
$$
\mathrm{rk} \widehat{HF} (S^3_{\alpha / \beta} (K)) =
|\alpha| + 2 \max( 0 , ( 2 \nu (K) - 1) |\beta| - |\alpha|) + |\beta| 
\left( \sum_s \left( \mathrm{rk} H_* ( \hat{A}_s ) - 1 \right) \right) .
$$
This holds for any pair of coprime integers $\alpha , \beta$. 

For simplicity, 
let $X(\nu (K), \alpha,\beta)$ denote $\max( 0 , ( 2 \nu (K) - 1) |\beta| - |\alpha|)$ and 
$Y$ denote $\sum_s \left( \mathrm{rk} H_* ( \hat{A}_s ) - 1 \right)$. 
Then we have 
\begin{equation}\label{rkHF}
\mathrm{rk} \widehat{HF} (S^3_{\alpha / \beta} (K)) =
|\alpha| + 2 X(\nu (K), \alpha,\beta) + |\beta| Y . 
\end{equation}

Now, for some coprime integers $\alpha , \beta$ with $\beta \ge 2$, 
we assume that $S^3_{\alpha / \beta} (K)$ is an L-space, 
i.e., by definition, 
$$\mathrm{rk} \widehat{HF} (S^3_{\alpha / \beta} (K)) = |\alpha|.$$
It then suffices to show that $S^3_{\alpha} (K)$ is an L-space, 
i.e., $\mathrm{rk} \widehat{HF} (S^3_{\alpha} (K)) = |\alpha|$.

On the other hand, in general, we see that 
$\mathrm{rk} \widehat{HF} (S^3_{\alpha} (K))  - |\alpha| \ge 0$ for any integer $\alpha$ as follows. 
In the proof of \cite[Proposition 5.1]{OzsvathSzabo1},
it is claimed that
$$\chi (\widehat{HF} (S^3_{\alpha} (K)) ) = | H_1 (S^3_{\alpha} (K) ;\mathbb{Z} )|. $$
Also see \cite[Section 2]{OzsvathSzabo2}.
By definition, the Euler characteristic (the left-hand side) is
the alternating sum of the dimensions of $\widehat{HF} (S^3_{\alpha} (K)) $. 
Hence, it is not greater 
than the total rank of $\widehat{HF} (S^3_{\alpha} (K)) $, i.e., 
$$\mathrm{rk} \widehat{HF} (S^3_{\alpha} (K))  \ge \chi (\widehat{HF} (S^3_{\alpha} (K)) ) 
= | H_1 (S^3_{\alpha} (K) ; \mathbb{Z} ) | = |\alpha|. $$

From this equation, 
in order to obtain $\mathrm{rk} \widehat{HF} (S^3_{\alpha} (K)) = |\alpha|$, 
it suffices to show that 
$\mathrm{rk} \widehat{HF} (S^3_{\alpha} (K))  - |\alpha| \le 0$. 
Actually, we have from equation~\eqref{rkHF}; 
\begin{equation}\label{rkHF1}
\mathrm{rk} \widehat{HF} (S^3_{\alpha} (K)) - |\alpha| = 2 X(\nu (K), \alpha, 1) + Y . 
\end{equation}

Note here that we have $Y \le 0$ as follows. 
It is seen that 
\begin{equation}\label{Lspace}
2 X(\nu (K), \alpha,\beta) + |\beta| Y =0 
\end{equation}
from equation~\eqref{rkHF} and 
the assumption that $\mathrm{rk} \widehat{HF} (S^3_{\alpha / \beta} (K)) = |\alpha|$. 
Thus, together with $X(\nu (K), \alpha,\beta) \ge 0$ by definition, we have $Y \le 0$.

If $\nu(K) \le 0$, then 
$$X(\nu (K), \alpha, 1) = \max( 0 , ( 2 \nu (K) - 1) - |\alpha|) =0 $$
holds. 
Since $Y \le 0$, together with equation~\eqref{rkHF1}, 
we obtain that $\mathrm{rk} \widehat{HF} (S^3_{\alpha} (K))  - |\alpha| \le 0$ as desired. 

If $\nu(K) \ge 1$, then we have 
$X(\nu (K), \alpha, 1) < X(\nu (K), \alpha, \beta)$ 
from the assumption that $\beta \ge2$. 
Thus, together with $Y \le 0$ and equation~\eqref{Lspace}, we obtain that
$$2 X(\nu (K), \alpha, 1) + Y < 2 X(\nu (K), \alpha, \beta) + Y = - |\beta| Y + Y \le 0 .$$
Together with equation~\eqref{rkHF1}, 
this implies that $\mathrm{rk} \widehat{HF} (S^3_{\alpha} (K))  - |\alpha| \le 0$ as desired. 

\end{proof}

Then, in \cite[Corollary 1.3]{OzsvathSzabo2}, 
Ozsv\'ath and Szab\'o proved that 
if a knot $K$ in $S^3$ admits an integral Dehn surgery 
yielding an L-space, 
then the Alexander polynomial $\Delta_K (t )$ has the form 
$$\Delta_K (t ) = (-1)^k + \sum^k_{j =1 } ( -1)^{k -j} 
\;\left(\;t^{n_j} + t^{- n_j} \;\right)$$
for some increasing sequence of positive integers 
$0 < n_1 < n_2 < \cdots < n_k$. 
This means that 
all non-zero coefficients of $\Delta_K(t)$ are $\pm 1$. 
\end{proof}

\begin{remark}
In the above proof, Claim~\ref{clm:Lsurgery} is actually necessary for the following reason. 
By the Cyclic Surgery Theorem established in \cite{CGLS}, 
all cyclic surgeries on hyperbolic knots in $S^3$ 
are shown to be integral surgeries. 
However, the Finite Surgery Theorem of \cite{BoyerZhang} 
shows that finite surgeries on hyperbolic knots in $S^3$ 
are half-integral or integral. 
In other words, at present, we cannot rule out the possibility of 
a half-integral surgery and it is currently only a conjecture that 
such finite surgeries are integral: 
See \cite[Problem 1.77 A(6)]{Kirby} for more information. 
\end{remark}

\section{Calculation of the Alexander polynomials}\label{sec3}

In this section, we prove the following proposition, 
which will be shown by direct calculations of the Alexander polynomials. 

\begin{proposition}\label{prop2}
Let $K$ be a pretzel knot of type $(-1, 2n, p, q)$, 
where $n$ is an integer and 
$p,q$ are odd positive integers with $p \le q$. 
If every non-zero coefficient of the Alexander polynomial of $K$
is $\pm 1$, then $n = 1$ and $p=3$. 
\end{proposition}

\begin{figure}[htb]
\includegraphics[width=120pt]{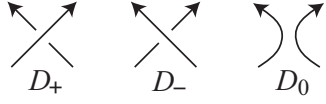}
\caption{Skein triples}\label{fig:skein}
\end{figure}

Recall that the Alexander polynomial $\Delta_L(t)$ for a link $L$ satisfies the following skein relation (see \cite[pp.~82]{Lickorish} for example): 
\begin{align}\label{eq:skein}
\Delta_{L_+}(t) - \Delta_{L_-}(t) = (t^{-1/2} - t^{1/2})\Delta_{L_0}(t),
\end{align}
where $L_+$, $L_-$, and $L_0$ possess diagrams $D_+$, $D_-$ and $D_0$ which differ only in a small neighborhood as shown in Figure~\ref{fig:skein}. 

\begin{remark}\label{rem:nabla}
Let $l$ be a positive integer, and $\Delta_l(t)$ the Alexander polynomial of a $(2,l)$-torus link. 
Set $f_l = \sum_{i=0}^{l}t^i$. 
Then we have $\Delta_l(-t) = (-t)^{(1-l)/2} f_{l-1}$. 
See \cite[pp.~98]{Kawauchi} for example. 
\end{remark}

\begin{proof}[Proof of Proposition \ref{prop2}]

We divide our proof of Proposition~\ref{prop2} into three claims. 
We denote by $P(a_1,\ldots,a_j)$ a pretzel link of type $(a_1,\ldots,a_j)$,
and by $[g(t)]_j$ the coefficient of $t^j$ in a polynomial $g(t)$. 

\bigskip

\begin{claim}\label{clm:-2n}
Let $n$ be an integer with $n \geq 1$. 
Let $p$ and $q$ be odd integers with $3 \leq p \leq q$. 
Let $K$ be a pretzel knot of type $(-1, -2n, p, q)$. 
Then we have 
\begin{align*}
[\Delta_K(t)]_1 =\left\{ \begin{array}{ll}
-4 & \text{if }~n=1, \\
-3 & \text{if }~n \geq 2, \\
\end{array} \right. 
\end{align*}
where $\Delta_K(t)$ is normalized so that $\mindeg \Delta_K(t) = 0$ and $[\Delta_K(t)]_0 > 0$.
\end{claim}
\begin{proof}
Let $K=P(-1, -2n, p, q)$ with $1 \leq n$ and $3 \leq p \leq q$. 
By applying the skein formula~\eqref{eq:skein} at crossings in the $(-2n)$-twists repeatedly, 
we can obtain a resolving tree such that each leaf node corresponds to either $P(-1, 0,p,q)$ or $P(-1, -1, p, q)$. 
Notice that $P(-1, 0, p, q)$ is equivalent to a connected sum of a $(2,p)$-torus knot and a $(2,q)$-torus knot. 
Then we have
\begin{align*}
\Delta_K(t) &= \Delta_{2n-1}(t) \Delta_{P(-1, 0, p, q)}(t) - \Delta_{2n}(t) \Delta_{P(-1, -1, p, q)}(t)\\
&= \Delta_{2n-1}(t) \Delta_{p}(t) \Delta_{q}(t) -  \Delta_{2n}(t) \Delta_{P(-1, -1, p, q)}(t).
\end{align*}
Next we calculate the Alexander polynomial of $P(-1,-1,p,q)$ by the same argument as above. 
By applying the skein formula~\eqref{eq:skein} at crossings in the $p$-twists repeatedly, 
we can obtain a resolving tree such that each leaf node corresponds to either $P(-1, -1,0,q)$ or $P(-1, -1, 1, q)$. 
Note that $P (-1,-1, 0, q)$ is equivalent to a $(2, q)$-torus knot and that $P(-1,-1,1,q)$ is equivalent to a $(2,q-1)$-torus link. 
Then we have
\begin{align*}
\Delta_{P(-1, -1,p,q)}(t) &= \Delta_{p-1}(t) \Delta_{P(-1,-1,0,q)}(t) + \Delta_{p}(t) \Delta_{P(-1,-1,1,q)}(t)\\
&= \Delta_{p-1}(t)\Delta_{q}(t) + \Delta_{p}(t) \Delta_{q-1}(t).
\end{align*}
Hence we have
\begin{align*}
\Delta_K(t) = \Delta_{2n-1}(t) \Delta_p(t) \Delta_q(t) - \Delta_{2n}(t) \Delta_{p-1}(t) \Delta_q(t) - \Delta_{2n}(t) \Delta_{p}(t) \Delta_{q-1}(t). 
\end{align*}
To calculate easily, we consider the polynomial obtained by substituting $-t$ in the Alexander polynomial, namely, $\Delta_K(-t)$. 
By Remark~\ref{rem:nabla}, we have 
\begin{align*}
\Delta_K(-t) &=(-t)^{(4-p-q-2n)/2} \left( f_{2n-2} f_{p-1} f_{q-1} - f_{2n-1} f_{p-2} f_{q-1} - f_{2n-1} f_{p-1} f_{q-2} \right)\\
&\doteq -f_{2n-2} f_{p-1} f_{q-1} + f_{2n-1} f_{p-2} f_{q-1} + f_{2n-1} f_{p-1} f_{q-2} .
\end{align*}
Here the symbol $\doteq$ means that both sides are equal up to multiplications by units of the Laurent polynomial ring $\mathbb{Z}[t,t^{-1}]$. 
Here we recall that $1 \leq n$ and $3 \leq p \leq q$. 
Then we have 
\begin{align*}
[f_{2n-2}f_{p-1}f_{q-1}]_1 &= \left\{ \begin{array}{ll}
2 & \text{if }~n=1, \\
3 & \text{if }~n \geq 2,\\
\end{array} \right. \\
[f_{2n-1}f_{p-2}f_{q-1}]_1 &= 3,\\
[f_{2n-1}f_{p-1}f_{q-2}]_1 &= 3.
\end{align*} 
Therefore we have 
\begin{align*}
[\Delta_K(-t)]_1 =\left\{ \begin{array}{ll}
-2+3+3 = 4 & \text{if }~n=1, \\
-3+3+3 = 3 & \text{if }~n \geq 2,\\
\end{array} \right. 
\end{align*}
that is, 
\begin{align*}
[\Delta_K(t)]_1 =\left\{ \begin{array}{ll}
-4 & \text{if }~n=1, \\
-3 & \text{if }~n \geq 2.\\
\end{array} \right. 
\end{align*}
\end{proof}

\begin{claim}\label{clm:+2n}
Let $n$ be an integer with $n \geq 2$.
Let $p$ and $q$ be odd integers with $3 \leq p \leq q$. 
Let $K$ be a pretzel knot of type $(-1,2n, p, q)$. 
Then we have 
\begin{align*}
[\Delta_K(t)]_3 = 2, 
\end{align*}
where $\Delta_K(t)$ is normalized so that $\mindeg \Delta_K(t) = 0$ and $[\Delta_K(t)]_0 > 0$.
\end{claim}
\begin{proof}
The proof is similar to that of Claim~\ref{clm:-2n}. 
Let $K=P(-1,2n, p, q)$ with $2 \leq n$ and $3 \leq p \leq q$. 
By applying the skein formula~\eqref{eq:skein} at crossings in the $2n$-twists repeatedly, 
we can obtain a resolving tree such that each leaf node corresponds to $P(-1,0,p,q)$ or $P(-1,1,p,q)$. 
Then we have
\begin{align*}
\Delta_K(t) &= \Delta_{2n-1}(t) \Delta_{P(-1,0,p,q)}(t) + \Delta_{2n}(t) \Delta_{P(-1,1,p,q)}(t)\\
&= \Delta_{2n-1}(t) \Delta_{p}(t) \Delta_{q}(t) + \Delta_{2n}(t) \Delta_{P(-1,1,p,q)}(t).
\end{align*}
By applying the same argument as above at crossings in the $p$-twists, we have
\begin{align*}
\Delta_{P(-1,1,p,q)}(t) &= \Delta_{p-1}(t) \Delta_{P(-1,1,0,q)}(t) + \Delta_{p}(t) \Delta_{P(-1,1,1,q)}(t)\\
&= \Delta_{p-1}(t)\Delta_{q}(t) + \Delta_{p}(t) \Delta_{P(-1,1,1,q)}(t).
\end{align*}
Notice that $P(-1,1,1,q)$ is equivalent to a $(2,q+1)$-torus link.
By applying the skein formula \eqref{eq:skein}, we have
$\Delta_{q+1}= \Delta_{q-1}(t) + (t^{-1/2}-t^{1/2})\Delta_{q}(t)$.
Hence we have
\begin{align*}
\Delta_K(t) &= \Delta_{2n-1}(t) \Delta_p(t) \Delta_q(t)  + \Delta_{2n}(t) \Delta_{p-1}(t) \Delta_q(t) + \Delta_{2n}(t) \Delta_{p}(t) \Delta_{q-1}(t)\\ &\hspace{1em}+ (t^{-1/2}-t^{1/2})\Delta_{2n}(t) \Delta_{p}(t) \Delta_{q}(t),
\end{align*}
and then we have
\begin{align*}
\Delta_K(-t) &\doteq -t f_{2n-2} f_{p-1}f_{q-1}- t f_{2n-1} f_{p-2}f_{q-1} - t f_{2n-1} f_{p-1}f_{q-2}\\ &\hspace{1em}+ (1 + t) f_{2n-1} f_{p-1}f_{q-1}. 
\end{align*}
Here we recall that $2 \leq n$ and $3 \leq p \leq q$. 
Then we have 
\begin{align*}
[t f_{2n-2} f_{p-1} f_{q-1}]_3 &= [f_{2n-2} f_{p-1} f_{q-1}]_2\\
&=6,\\
[t f_{2n-1} f_{p-2}f_{q-1}]_3 &= [f_{2n-1} f_{p-2} f_{q-1}]_2\\
&= \left\{ \begin{array}{ll}
5 & \text{if }~p=3,~ q \geq 3,\\
6 & \text{if }~5 \leq p \leq q,
\end{array} \right.\\
[t f_{2n-1} f_{p-1} f_{q-2}]_3 &= [f_{2n-1} f_{p-1} f_{q-2}]_2\\
&= \left\{ \begin{array}{ll}
5 & \text{if }~p=3,~q=3, \\
6 & \text{if }~p \geq 3,~q \geq 5,
\end{array} \right.\\
[(1+t) f_{2n-1} f_{p-1} f_{q-1}]_3 &= [f_{2n-1} f_{p-1} f_{q-1}]_3 + [f_{2n-1} f_{p-1} f_{q-1}]_2\\
&= \left\{ \begin{array}{ll}
8+6 =14 & \text{if }~p=3,~q=3, \\
9+6 =15 & \text{if }~p=3,~q \geq 5,\\
10+6 =16 & \text{if }~5 \leq p \leq q.
\end{array} \right.
\end{align*}
Therefore we have 
\begin{align*}
[\Delta_K(-t)]_3&= \left\{ \begin{array}{ll}
-6 -5 -5 +14 =-2 & \text{if }~p=3,~q=3, \\
-6 -5 -6 +15 =-2 & \text{if }~p=3,~q \geq 5,\\
-6 -6 -6 +16 =-2 & \text{if }~5 \leq p \leq q,
\end{array} \right.
\end{align*}
that is, $[\Delta_K(t)]_3=2$. 
\end{proof}

Here we note that $P(-1, 2, p, q)$ is equivalent to $P(-2, p, q)$.

\begin{claim}\label{clm:-2pq}
Let $p$ and $q$ be odd integers with $5 \leq p \leq q $. 
Let $K$ be a pretzel knot of type $(-2, p, q)$. 
Then we have 
\begin{align*}
[\Delta_K(t)]_4 = -2, 
\end{align*}
where $\Delta_K(t)$ is normalized so that $\mindeg \Delta_K(t) = 0$ and $[\Delta_K(t)]_0 > 0$.
\end{claim}
\begin{proof}
Let $K=P(-2,p,q)$ with $5 \leq p \leq q$. 
By applying the skein formula~\eqref{eq:skein} at a crossing in the $(-2)$-twists, 
we have
\begin{align*}
\Delta_K(t) = \Delta_p(t) \Delta_q(t) + (t^{-1/2} - t^{1/2}) \Delta_{p+q}(t).
\end{align*}
Then we have
\begin{align*}
\Delta_K(-t) \doteq -tf_{p-1}f_{q-1} +(1+t)f_{p+q-1}. 
\end{align*}
Here we recall that $5 \leq p \leq q$. 
Then we have 
$[tf_{p-1}f_{q-1}]_4 = [f_{p-1}f_{q-1}]_3 =4$ and $[(1+t)f_{p+q-1}]_4 = 1 + 1 = 2$. 
Therefore we have $[\Delta_K(-t)]_4 = -4 +2 = -2$, 
that is, $[\Delta_K(t)]_4 = -2$.
\end{proof}

This completes the proof of Proposition \ref{prop2}.
\end{proof}

\section{Proof of Theorem \ref{MainThm}}\label{sec4}

\begin{proof}[Proof of Theorem \ref{MainThm}]
By Propositions \ref{prop1} and \ref{prop2}, 
if a hyperbolic Montesinos knot $K$ 
admits a non-trivial cyclic or finite surgery, 
then $K$ is equivalent to a $(-1, 2, 3, q)$-pretzel knot, 
where $q$ is an odd positive integer with $3 \leq q$. 
This $K$ is actually equivalent to a $(-2, 3, q)$-pretzel knot. 
Then Mattman showed in \cite[Theorem 1.1 and 1.2]{Mattman} that, 
among such pretzel knots, 
only the $(-2, 3, 7)$- and $(-2, 3, 9)$- can have cyclic/finite surgeries, 
and the surgery slopes are as described in Theorem \ref{MainThm}. 
This completes the proof of Theorem \ref{MainThm}.
\end{proof}

\begin{remark}
The techniques we have used in this paper cannot be applied to 
the $( -2, 3, q)$-pretzel knots as they are fibered 
and all non-zero coefficients of their Alexander polynomials are $\pm 1$.
\end{remark}


\bibliographystyle{amsplain}

\end{document}